\documentclass[%
reprint,
 amsmath,amssymb,
 aps,
]{revtex4-2}
\usepackage{float, xcolor}
\usepackage{wrapfig}
\usepackage{graphicx}
\usepackage{dcolumn}
\usepackage{soul}
\usepackage{cancel}
\usepackage{bm}
\usepackage{subfigure}
\usepackage{ulem} 
\usepackage{threeparttable} 
\usepackage{wrapfig}
\usepackage{graphicx}
\usepackage{dcolumn}
\usepackage{ulem}
\usepackage[colorlinks,citecolor=blue,urlcolor=blue,linkcolor=blue,breaklinks=true]{hyperref}
\usepackage[mathlines]{lineno}
\usepackage[marginal]{footmisc}

\begin{document}

\preprint{APS/123-QED}

\title{$\alpha$-SGHN: A Robust Model for Learning Particle Interactions in Lattice Systems
}

\author{Yixian Gao \textsuperscript{1} }
\email{These authors contributed equally to this work.}

\author{Ru Geng\textsuperscript{1}}
\email{These authors contributed equally to this work.}

\author{Panayotis Kevrekidis\textsuperscript{2}}

 \author{Hong-Kun Zhang\textsuperscript{2}}
\email{corresponding author, hongkun@math.umass.edu}

\author{Jian Zu \textsuperscript{1}}
\email{corresponding author, zuj100@nenu.edu.cn}

\address{%
1.Center for Mathematics and Interdisciplinary Sciences,
School of Mathematics and Statistics,  Northeast Normal University,  Changchun, 130024, P.R. China
}%
\address{%
2.Department of Mathematics and Statistics, University of Massachusetts,   Amherst, 01003, MA, USA
}%

\date{\today}

\begin{abstract}
 
 We propose an $\alpha$-separable graph Hamiltonian network ($\alpha$-SGHN) that reveals complex interaction patterns between particles in lattice systems. Utilizing trajectory data, $\alpha$-SGHN infers potential interactions without prior knowledge about particle coupling, overcoming the limitations of traditional graph neural networks that require predefined links. Furthermore, $\alpha$-SGHN preserves all conservation laws during trajectory prediction. Experimental results demonstrate that our model, incorporating structural information, outperforms baseline models based on conventional neural networks in predicting lattice systems. We anticipate that the results presented will be applicable beyond the specific onsite and inter-site interaction lattices studied, including the Frenkel-Kontorova model, the rotator lattice, and the Toda lattice.

 \end{abstract}

\maketitle

\begin{figure*}  %
\centering
\includegraphics[scale=0.5]{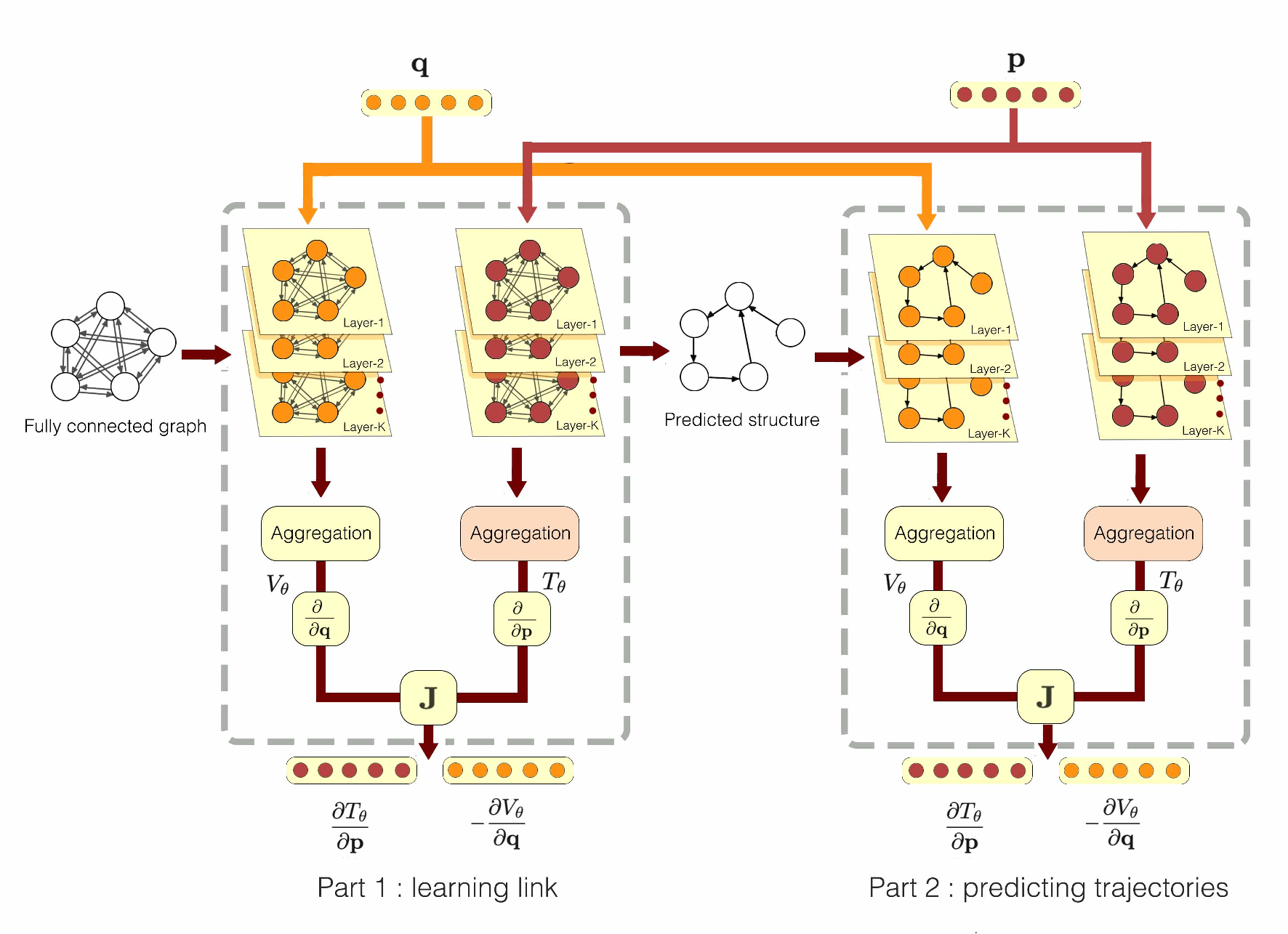}
\caption{ The architecture of $\alpha$-SGHN}\label{alphaSGHN}
\end{figure*} 
\section{Introduction}\label{sec1}

\begin{table*}
\caption{\label{com_net}%
Network performance comparison. $\alpha$-SGHN has the least input conditions and the most widely applicable prediction content.
}
\begin{ruledtabular}
\begin{tabular}{c|c|c|c}
\textrm{}&\multicolumn{1}{c|}{\textrm{Additional information required}}&\multicolumn{1}{c|}{\textrm{Predicting particle interaction }} &\multicolumn{1}{c|} System potential energy requirements 
\\
\colrule
MLP &   -  &  no &  any \\
HNN  &  -  &  no &  any   \\
HOGN &  Particle interaction relationship &  no  &  even symmetry  \\
HGNN & Particle interaction relationship &  no   &  even symmetry  \\
$\alpha$-SGHN(ours) &-  &  yes  & any  \\
 \end{tabular}
\end{ruledtabular}
\end{table*}

The impact of lattice systems on the discoveries and technological progress of condensed matter physics, materials science, biochemistry, medicine, and other scientific fields is profound.
 For example,  Fermi, Pasta, Ulam and 
 Tsingou attempted to study the thermalization process through a prototypical nonlinear lattice~\cite{fermi1955studies}.
The Frenkel–Kontorova system is commonly used to study heat conduction theory~\cite{hu2005heat}, DNA dynamics~\cite{braun2004frenkel,Yomosa1983,homma1984coupled,zhou1991short}, lattice defects/dislocations and crowdions \cite{braun2004frenkel,kovalev1994one,landau1993model},  hydrogen-bonded chains\cite{antonchenko1983solitons}, etc.; see
also~\cite{SGbook}.
The $\phi^4$ system describes possible domain walls in early universe models in cosmology \cite{zel1974cosmological}, studies structural phase transitions in the displacive limit in condensed matter physics \cite{schneider1975observation,currie1980statistical,mazenko1978statistical}, and serves as a phenomenological model for nonlinear excitations  with exotic spin-charge relations \cite{su1979solitons,jackiw1981solitons}, etc.;
a recent recap of the relevant model and
applications can be found in~\cite{p4book}.

These types of lattice systems can be generally described by the Hamiltonian
\begin{equation}\label{dpq}
\frac{d}{dt}\left(
\begin{array}{c} 
\mathbf{q}\\
\mathbf{p}
\end{array}
\right)=
\mathbf{J}\left(
\begin{array}{c} 
\frac{\partial {H}}{\partial \mathbf{q}}\\
\frac{\partial {H}}{\partial \mathbf{p}}
\end{array}
\right),\quad \mathbf{J}=\left(
\begin{array}{c c} 
O &  \mathbf{I}\\
- \mathbf{I} &O
\end{array}\right),
\end{equation}
where $\mathbf{q}=(q_1, \cdots, q_N)$ and $\mathbf{p}=(p_1, \cdots, p_N)$ are the generalized positions and momenta, respectively, of the system of point-like particles; $N$ is the number of  point-like particles in the system.
 $\mathbf{I}$ represents the identity matrix. 
Here  $H=T+V$ is the Hamiltonian, with  kinetic energy $$ {T}(\mathbf{p})=\sum_{i=1}^N T({p_i})=\sum_{i=1}^N \frac{p_{i}^2 }{2m}$$ and potential
\begin{eqnarray}\label{phieq}
V(\mathbf{q})
=\sum_{i=1}^{N} \left( V^{(1)}(q_{i+1}-q_i) +V^{(2)}(q_i) \right).
\end{eqnarray}
The inter-site potential $V^{(1)}$ encapsulates particle interactions within the nearest-neighbour, while the on-site potential $V^{(2)}$ accounts for potential interactions with an external environment, such as a substrate.
 Moreover  $m$ denotes the mass of the point-like particle, and we set $m=1$ for simplicity.  

Finding solutions to the differential equations based only on the observed trajectories is an important task in scientific fields.  The emergence of artificial neural networks represents a revolutionary development in the analysis and interpretation of complex data. Paired or multi particle energy interactions in lattice systems can lead to interesting phenomena in the system's behavior, such as phase transitions, phonon vibrations, magnetic behavior, soliton formation, and 
etc.~\cite{SGbook,p4book,kevrekidis2020emerging}. In the specific context of lattice systems, machine learning methods, especially deep learning, have made significant progress in describing governing equations, identifying phase transitions, and constructing predictive models \cite{raissi2019physics,lu2021deepxde,Kutz_SINDy,saqlain2022discovering, zvyagintseva2022machine, li2022gradient, zhu2022neural, jin2022learning}. Although they are useful, the methods that rely on conventional neural networks still pose challenges. Some require knowledge or understanding of system equations \cite{li2022gradient,zhu2022neural}, while others require creating an overcomplete operator library \cite{Kutz_SINDy,raissi2019physics,saqlain2022discovering}.

On the other hand, data that exhibits relationships between elements can be modeled as a graph. Elements are named graph nodes, and relationships are edges. A graph can naturally serve as a comprehensive representation of a lattice system, where the particles of the system correspond to the nodes of the graph, and the interactions between particles are represented as edges of the graph. Therefore, Graph Neural Networks (GNNs) have become a new field for studying lattice systems.
When dealing with lattice systems involving many degrees of freedom and complex interactions, compared to conventional neural networks, GNNs exhibit a proficient ability to accurately identify key information around nodes, which ultimately improves the accuracy of the model.  The outstanding performance of graph neural networks in some classic systems can be observed, for example, in the works cited in \cite{sanchez2019hamiltonian,bishnoi2023learning}.
An important advantage of graph neural networks is that they encode the interaction relationships of the system (i.e., the underlying structure of the system). 
In some  special cases, lattice systems such as gravitationally 
interacting celestial bodies,  can be modeled using fully connected graphs (due to their interactions existing between all particles).
The interactions in most lattice problems are complex and diverse, and may exhibit irregularities \cite{2010Automated} or even demonstrate long-range interactions \cite{laskin2006nonlinear,go1978respective}, and cannot be observed in advance.
This complexity makes visual recognition of which particles interact (i.e.,
between which a link exists) challenging. Therefore, although methods based on graph neural networks are effective, they require more input information, especially the linking of system particles, compared to traditional neural networks.

To address these challenges, based on the previous work of a subgroup of the present authors~\cite{geng2024separable}, we propose the  $\alpha$-Separable
Graph Hamiltonian Network ($\alpha$-SGHN) model, which can extract the underlying particle link (associated with each particle's interactions) of the lattice system based solely on its motion trajectory facilitating subsequent system predictions or other applications. To our knowledge, this is 
a prototypical example that can provide the underlying particle link of a lattice Hamiltonian
system solely based on its trajectories.
Table \ref{com_net} compares our method with established models such as  multi-layer perceptrons (MLP) \cite{goodfellow2016deep} and Hamiltonian Neural Networks (HNNs) \cite{greydanus2019hamiltonian}, Hamiltonian ODE graph network (HOGN)\cite{sanchez2019hamiltonian},
and Hamiltonian graph neural networks (HGNN)\cite{bishnoi2023learning}.  Among them, HOGN and HGNN are graph method based. Our method only requires trajectory data to provide the most widely applicable predictive content. See Section \ref{Experiment} for specific comparison.
Considering that PINN and symbolic regression methods require prior knowledge of the governing equations, they are not included in our comparison.

Our model is divided into two parts. The first part is to learn the underlying particle interactions (links) in the system, and the second part is to apply the underlying graph structure (including link information) for more accurate prediction, as shown in Fig. \ref{alphaSGHN}. In addition, we investigate  whether the particle behavior predicted by $\alpha$-SGHN  preserves the system's conservation laws. The experiment shows that due to the embedding of particle interaction information in the GNN, it outperforms the baseline model based on conventional neural networks. In the appendix, we also explored the predictive effect of $\alpha$-SGHN on long-range interactions and high dimension system structures.  

The main contributions of our article are as follows:

1. Under the assumption of only knowing the trajectory of the system, our model learned the underlying interactions of the system particles. This addressed the limitation of GNN applications requiring knowledge of graph structure information.

2. Our model maintains the system's conservation laws during trajectory prediction.

3. We provide an explanation about why we anticipate that graph neural network models will generally outperform traditional neural network models when dealing with lattice systems.

4. We provide a method for modeling directed graphs in learning even non-symmetric potential energy systems.

Our presentation is structured as follows.
In section II we present our methodology.
Our numerical experiments are analyzed in detail
in section III. In section IV we make a brief
excursion to the intriguing realm of complete
integrability (in lattices). Finally, in section V
we summarize our findings and present our conclusions,
as well as some interesting directions for further
work.

\section{Experiment} \label{Experiment}
Our method is based on the previous work of~\cite{geng2024separable}, but the difference is that in this work, we can learn the interaction relationship (edge link) between particles in the lattice system. Please refer to the supplementary materials for specific details.
In this section, we use three test examples: the Frenkel-Kontorova lattice (one conserved quantity), the rotator lattice (two conserved quantities), and the Toda lattice (multiple ---i.e., as many 
as the degrees of freedom--- conserved quantities).
To ascertain the effectiveness of our proposed model's prediction performance, we conduct experiments comparing $\alpha$-SGHN with established network structures such as the MLP and the HNN.
Due to our assumption that the particle linkage of the system is unknown, i.e. we don't know the graph structure of the system, the existing graph neural network models cannot serve as our primary baseline model as they all require knowledge of the exact particle link. But we also present the results of learning HOGN and HGNN using fully connected graphs (assuming that there are interactions between all particles) in Table \ref{FK_com}.
We will verify through experiments that our graph neural network model $\alpha$- SGHN can effectively learn multiple conserved quantities due to its ability to learn potential particle linkage relationships in the system, and is superior to conventional neural network models. 
\subsection{ Lattice system}
\subsubsection{Frenkel-Kontorova lattice}

The Frenkel-Kontorova (FK) model was introduced  by Ya. Frenkel and T. Kontorova as a prototypical way to describe the structure and dynamics of a crystalline lattice near a dislocation. 
In solid-state physics, the model is, arguably, one
of the first instances where a two-dimensional bulk defect is represented by a straightforward one-dimensional chain~\cite{tekic2016ac}. 
The characteristics of the FK model provide deep physical insights and significantly simplify the understanding of nonlinear dynamics in various problems of solid-state physics, classical mechanics
and biophysics among many others~\cite{braun2004frenkel,SGbook}.
Its  Hamiltonian is given by
\begin{eqnarray}
H=\sum_{i=1}^N \Bigg(\frac{p_i^2}{2}+ \frac{(q_{i+1}-q_i)^2}{2}+1-\cos(q_i)\Bigg).	
\end{eqnarray}

\subsubsection{Rotator lattice}

The rotator lattice is, arguably, one of the simplest examples of a classical spin 1D model with nearest neighbor interactions, and has a potential function governed by $1-\cos(q_{i+1}-q_i)$.  This model can also be read as an array of $N$-coupled pendulum. 
As an example of a chaotic dynamical system, when it is reduced to a harmonic chain and free rotator, it is integrable in both  small and high energy  limits and has been widely studied~\cite{Giardina2000}.

The Hamiltonian for the coupled rotator lattice dynamics amended with harmonic interactions 
reads~\cite{Li_2015}:

\begin{eqnarray}
H=\sum_{i=1}^N \Bigg(\frac{p_i^2}{2}+ \frac{(q_{i+1}-q_i)^2}{2}+1-\cos(q_{i+1}-q_i)\Bigg).	
\end{eqnarray}

It has two conserved quantities of energy and momentum, with the latter being equal
to: $$P=\sum_{i=1}^N p_i =\sum_{i=1}^N \dot{q}_i$$.

 \subsubsection{Toda lattice}

The Toda lattice model~\cite{toda1967,toda1969,toda1970,Toda} has garnered significant attention in the physical sciences as a paradigm of a nonlinear lattice system where specific (solitonic) waveforms are known to propagate unaltered in shape, a property rigorously demonstrated in the literature \cite{ford1973integrability}. The system's dynamics are governed by the Hamiltonian:
\begin{eqnarray}\label{Toda}
H=\sum_{i=1}^N \Bigg(\frac{p_i^2}{2}+\exp(q_i-q_{i+1})\Bigg).
\end{eqnarray}
This system is known to feature $N$ constants of motion~\cite{flaschka1974toda}.
when the boundary conditions (BCs) are periodic, i.e. $q_{i+N}=q_i$ and $p_{i+N}=p_i$. In this case,
the $2 N$ phase-space coordinates can be used to define a symmetric, periodic tridiagonal $N \times N$ (time dependent) matrix,
$$
L=\left(\begin{array}{ccccc}
p_1 & v_1 & 0 & \cdots & v_{N} \\
v_1 & p_2 & v_2 & \cdots & 0 \\
0 & v_2 & p_3 & \cdots & 0 \\
\vdots & \vdots & \ddots & \ddots & \vdots \\
v_{N} & 0 & \cdots & v_{N-1} & p_N
\end{array}\right),
$$
where $v_i=-\mathrm{e}^{\left(q_a-q_{a+1}\right) / 2}$. All the quantities
$$
\mathcal{C}_n \equiv \sum_{i=1}^{\mathcal{N}} \lambda_i^n=\operatorname{Tr} L^n,
$$
then represent the constants of the motion, see \cite{flaschka1974toda} for details. 
The first two are familiar and physically meaningful, namely they represent the system's momentum  and energy.

\subsection{Dataset acquisition}

For FK and rotator lattice, we set $N=32$. For the Toda lattice, we set $N=3, 4$, and $5$ to verify the conservation laws of the system in the
cases of 3, 4, and 5 conserved quantities, respectively.
We use the explicit, symplectic Runge-Kutta-Nystr\"om algorithms with five stages (Section 8.5.3 in \cite{sanz2018numerical}) and time step 0.0025  to find 50 trajectories.
The initial conditions (ICs) are
 \begin{eqnarray}
&&\textbf{q}_0(i)\sim \mathcal{U}(0,1),\quad i=1,\cdots,N,\\
&&\textbf{p}_0(i)\sim \mathcal{U}(0,1),\quad i=1,\cdots,N,
\end{eqnarray}
where $ \mathcal{U}$ represents the uniform distribution.  Their BCs are periodic, i.e. $q_{i+N}=q_i$ and $p_{i+N}=p_i$.
We sub-sample the trajectories at a fixed time-step of 0.05 as the training set.

We use the same method to generate 20 trajectories as the test set, where the integration time is three times that of the training set, i.e. 15 time units.

\subsection{Network model settings}
For the baseline model MLP and HNN, we carefully select the hyperparameters of the model to achieve the best predictive performance. We start from a layer of 200 hidden units, continuously increasing and decreasing the network width and depth, and replacing different activation functions tanh, silu, and gelu until the network performance reaches 
what we estimated as an optimal level.
We carefully select the parameters of the baseline model to achieve the best performance.  We will demonstrate through experiments that $\alpha$-SGHN is far superior to the baseline model's optimal performance.

We adopt a learning rate piecewise constant decay strategy \cite{montavon2012neural},
see Table \ref{learning_rate}.
The piecewise learning rates are divided by 3500 and  5000 epoch point.
We record the loss function during the training process to ensure that the loss function of the neural network model reached steady convergence during the training process.
The total epoch for all model training is set to 10000.
The optimizer is Adam, and the batch size is 256.

\begin{table*}[h]
  \caption{Piecewise learning rate.}
  \label{learning_rate}
  \centering
  \begin{ruledtabular}
  \begin{tabular}{|c|c|c|c|c|c|}
    &FK         & Rotator  & Toda-3 & Toda-4  & Toda-5\\
\colrule
MLP-gelu       &    $10^{-3}$, $10^{-4}$, $10^{-5}$      &$10^{-3}$, $10^{-4}$, $10^{-5}$   &$10^{-3}$, $10^{-4}$, $10^{-5}$  &$10^{-3}$, $10^{-4}$, $10^{-5}$   & $10^{-3}$, $10^{-4}$, $10^{-5}$  \\
MLP-silu       &    $10^{-2}$, $5\times 10^{-3}$, $10^{-4}$      &$10^{-3}$, $10^{-4}$, $10^{-5}$   & $10^{-2}$, $5\times 10^{-3}$, $10^{-4}$  & $10^{-2}$, $5\times 10^{-3}$, $10^{-4}$  &   $10^{-2}$, $5\times 10^{-3}$, $10^{-4}$   \\
MLP-tanh      &    $10^{-3}$, $10^{-4}$, $10^{-5}$      &$10^{-3}$, $10^{-4}$, $10^{-5}$   &$10^{-3}$, $10^{-4}$, $10^{-5}$  &$10^{-3}$, $10^{-4}$, $10^{-5}$   &  $10^{-3}$, $10^{-4}$, $10^{-5}$   \\
\colrule
HNN-gelu   &    $10^{-3}$, $10^{-4}$, $10^{-5}$      &$10^{-3}$, $10^{-4}$, $10^{-5}$   &$10^{-3}$, $10^{-4}$, $10^{-5}$  &$10^{-3}$, $10^{-4}$, $10^{-5}$   &  $10^{-3}$, $10^{-4}$, $10^{-5}$   \\
HNN-silu   &   $10^{-2}$, $5\times 10^{-3}$, $10^{-4}$      &$10^{-3}$, $10^{-4}$, $10^{-5}$   & $10^{-2}$, $5\times 10^{-3}$, $10^{-4}$  & $10^{-2}$, $5\times 10^{-3}$, $10^{-4}$   &   $10^{-2}$, $5\times 10^{-3}$, $10^{-4}$   \\
HNN-tanh   &     $10^{-3}$, $10^{-4}$, $10^{-5}$      &$10^{-3}$, $10^{-4}$, $10^{-5}$   &$10^{-3}$, $10^{-4}$, $10^{-5}$  &$10^{-3}$, $10^{-4}$, $10^{-5}$   &  $10^{-3}$, $10^{-4}$, $10^{-5}$   \\
\colrule
$\alpha$-SGHN   &  $10^{-3}$, $10^{-4}$, $10^{-5}$      &$10^{-3}$, $10^{-4}$, $10^{-5}$   &$10^{-3}$, $10^{-4}$, $10^{-5}$  &$10^{-3}$, $10^{-4}$, $10^{-5}$   &  $10^{-3}$, $10^{-4}$, $10^{-5}$   \\
\end{tabular}
\end{ruledtabular}
\end{table*}

\subsection{Test metrics}

To evaluate our model in the test set, we log the following metrics:
Testing loss, the Mean Square Error (MSE) of the predicted trajectories, the MSE of the predicted conserved quantities featured in the system.

To determine the predicted conserved quantities and trajectories metric over long timespans, 
we integrate neural network models according to \eqref{pred}
by explicit, symplectic Runge-Kutta-Nystr\"om algorithms with time step 0.0025. $(\mathbf{q}^0,\mathbf{p}^0)$ is the initial value in the test set. In our settings, the initial time $t_0=0s$. 
While we are cognizant of the interesting recent work
of \cite{danieli2024dynamical} which suggests
that even such symplectic integrators can yield
adverse (integrability or conservation law-breaking
features), this is only relevant for times much longer
(and time steps much bigger) than considered herein.
 
The MSE of the predicted trajectories is defined as
\begin{eqnarray*}
\rm{MSE}_{traj}=\sum_{i=1}^N \bigg((q_{i}^t-{\hat{q}}_{i}^t)^2+(p_{i}^t-{\hat{p}}_{i}^t)^2\bigg).
\end{eqnarray*}

 \subsection{Link Extraction}\label{Structural_extraction}
 
\begin{figure*}%
\centering
\includegraphics[scale=0.6]{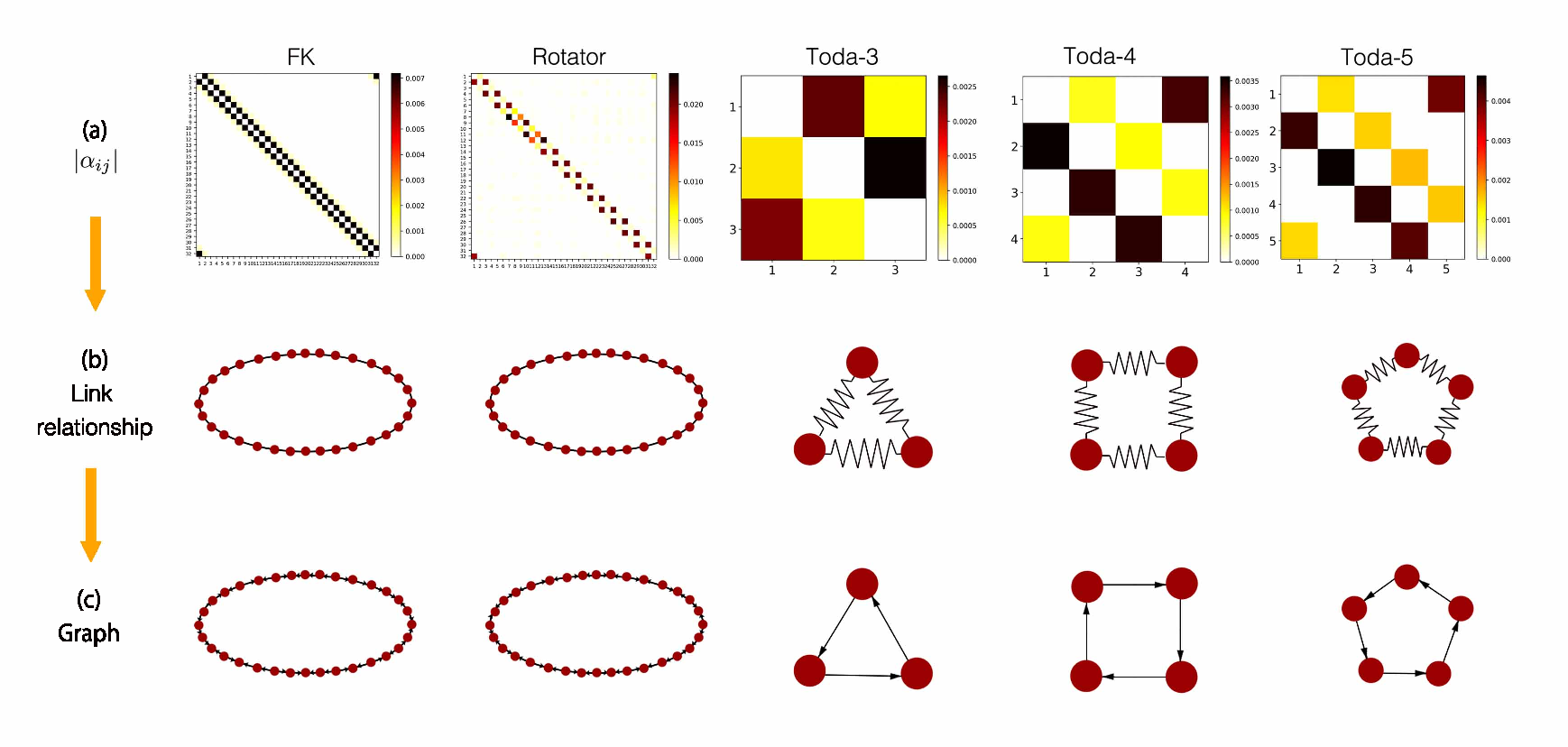}
\caption{(a). The $x$-axis and $y$-axis represent nodes, and the color coded $|\alpha_{i,j}|$ values. (b) is the particle link relationship extracted from (a). Obviously, it satisfies the periodic boundary conditions. (c) is a directed graph constructed through (b), with directions that are relatively large (or small) between $|\alpha_{i,j}|$ and $|\alpha_{j,i}|$.}\label{fig_edges}
\end{figure*}

 Fig. \ref{fig_edges} (a) illustrates $|\alpha_{i,j}|$, where the $x$-axis and $y$-axis represent nodes, and the color indicates the strength of the interact between nodes. We can extract the linking relationship of system particles from  Fig. \ref{fig_edges} (a), as shown in  Fig. \ref{fig_edges} (b). Based on the linkage relationship, we further construct a unidirectional graph $\mathcal{G}_{\alpha}$, where the direction of the graph edges goes towards the direction with higher (or lower) $|\alpha_{i,j}|$ values, as shown in  Fig. \ref{fig_edges} (c).

\subsection{Prediction Results}
Table \Ref {FK_com} compares the predictive performance of $\alpha$-SGHN and baseline models for FK system. We trained the baseline model to its optimal performance and recorded the effects of changing the depth, width, and activation function of the network, which did not further improve the performance of the baseline model.
It is evident that $\alpha$-SGHN outperforms the baseline models considered herein in terms of testing loss, long-term trajectory prediction, and energy conservation properties.
This is because the interior of the lattice system largely depends on the interactions between particles.
Conventional neural networks can only process data for each particle equivalently, making it difficult to find interactions between particles. $\alpha$-SGHN can effectively predict the particle interactions within the system and learn the lattice system more effectively by utilizing the predicted interactions.

\begin{table*}
\caption{\label{FK_com}%
The test results for FK system with $N=32$.  The best results are emphasized by bold fonts.  
}
\begin{ruledtabular}
\begin{tabular}{cccc}
\textrm{}&
\textrm{Test loss}&
\textrm{Trajectory$_{MSE}$}&
\textrm{Energy$_{MSE}$}\\
\colrule
MLP-1-100-silu   &  3.70E-3 $\pm$ 1.20E-3   &   1.32E-1 $\pm$ 1.43E-1  &  1.56E+0 $\pm$  1.83E+0 \\
MLP-1-100-tanh  &  2.36E-3 $\pm$ 8.52E-4   &   8.81E-2 $\pm$ 1.04E-1  &  3.70E-1 $\pm$  3.55E-1 \\
MLP-1-50-gelu &  3.31E-2 $\pm$ 5.95E-3   &   7.91E-1 $\pm$ 8.00E-1  &  9.57E+1 $\pm$  1.90E+2 \\
MLP-2-100-gelu &  3.58E-3 $\pm$ 9.57E-4   &   7.30E-2 $\pm$ 7.45E-2&  4.18E-1 $\pm$  3.61E-1 \\
MLP-1-200-gelu  &  1.89E-2 $\pm$ 4.72E-3   &   6.15E-1 $\pm$ 7.14E-1&  1.70E+2 $\pm$  2.79E+2 \\
MLP-1-100-gelu  & \uwave{2.06E-3 $\pm$ 6.09E-4 }  &    \uwave{3.91E-2 $\pm$ 3.72E-2} & \uwave{1.33E-1 $\pm$  1.17E-1} \\
\colrule
HNN-1-100-silu &  1.91E-3 $\pm$ 5.73E-4   &   4.57E-2 $\pm$ 4.27E-2  &  6.76E-3 $\pm$  4.64E-3 \\
HNN-1-100-tanh   &  1.66E-3 $\pm$ 6.79E-4   &   5.06E-2 $\pm$ 5.40E-2  &  9.34E-3 $\pm$  1.21E-2 \\
HNN-1-50-gelu &  3.68E-2 $\pm$ 6.57E-3   &   7.62E-1 $\pm$ 9.03E-1  &  2.53E+2 $\pm$  5.30E+2 \\
HNN-2-100-gelu  &  9.06E-4 $\pm$ 2.90E-4   &   1.41E-2 $\pm$ 1.31E-2  &  3.29E-3 $\pm$  2.90E-3 \\
HNN-1-200-gelu  &  1.43E-3 $\pm$ 3.99E-4   &   3.20E-2 $\pm$ 3.08E-2  &  4.61E-3 $\pm$  3.36E-3 \\
HNN-1-100-gelu & \uwave{7.44E-4 $\pm$ 2.62E-4 }  &   \uwave{1.36E-2 $\pm$ 1.29E-2}  & \uwave{2.21E-3 $\pm$  1.65E-3} \\
\colrule
 $\alpha$-SGHN(ours) &\textbf{4.26E-7} $\pm$ \textbf{1.22E-6}  &\textbf{4.85E-6} $\pm$ \textbf{2.37E-5}  & \textbf{9.70E-7} $\pm$ \textbf{3.13E-6}\\
\end{tabular}
\end{ruledtabular}
\end{table*}

 It can also be seen from Table \ref{Rotator_com}  that the performance of $\alpha$-SGHN is still far superior to the baseline models we examined
in the case of the rotator lattice.
In this case, there exist two conservation laws
and both of them are retrieved by numerous
orders of magnitude better in the present
setting rather than by the MLP or the HNN.
\begin{table*}
\caption{\label{Rotator_com}%
The test results for Rotator system with $N=32$.  The best results are emphasized by bold fonts. }
\begin{ruledtabular}
\begin{tabular}{ccccc}
\textrm{}&
\textrm{Test loss}&
\textrm{Trajectory$_{MSE}$}&
\textrm{Energy$_{MSE}$}&
\textrm{Momentum$_{MSE}$}\\
\colrule
MLP-1-100-gelu  &  1.00E-3 $\pm$ 3.94E-4   &   1.93E-2 $\pm$ 3.20E-2  &  7.00E-2 $\pm$  7.27E-2  &  1.03E-2 $\pm$ 1.03E-2\\
MLP-1-100-tanh  &  9.90E-4 $\pm$ 4.67E-4   &   4.30E-2 $\pm$ 6.24E-2  &  5.22E-2 $\pm$  5.09E-2  &  \uwave{8.28E-6  $\pm$  9.70E-6}\\
MLP-1-50-silu  &  1.36E-2 $\pm$ 3.68E-3   &   6.01E-1 $\pm$ 5.82E-1 & 9.52E+1 $\pm$ 1.57E+2  & 8.06E-5 $\pm$ 1.41E-5\\
MLP-2-100-silu  &   1.31E-3 $\pm$ 4.91E-4    &  1.53E-2 $\pm$ 2.68E-2  &  7.49E-2 $\pm$ 8.44E-2   & 2.56E-2 $\pm$ 2.33E-2 \\
MLP-1-200-silu  &  1.11E-2 $\pm$ 2.32E-3   &  2.35E-1 $\pm$ 2.72E-1  &  1.76E+1 $\pm$ 1.73E+1   & 7.48E-5 $\pm$ 9.21E-5 \\
MLP-1-100-silu  &  \uwave{6.65E-4 $\pm$ 1.94E-4  }  &   \uwave{6.14E-3 $\pm$ 8.84E-3}  &   \uwave{3.15E-2 $\pm$ 2.86E-2 }  & 3.95E-2 $\pm$ 2.58E-2 \\
\colrule
\colrule
HNN-1-100-silu &  8.61E-4 $\pm$ 2.25E-4   &   \uwave{1.43E-2 $\pm$ 1.90E-2}  &  2.84E-2 $\pm$  1.46E-2 &  1.17E-1 $\pm$  3.65E-2\\
HNN-1-100-tanh &   1.47E-3 $\pm$ 6.57E-4  &    4.76E+0 $\pm$ 2.94E+1  &   9.50E+4 $\pm$  4.14E+5& 1.46E+1 $\pm$  6.35E+1\\
HNN-1-50-gelu  &   1.61E-2 $\pm$ 4.03E-3  &   2.63E+1 $\pm$ 1.60E+2 &   7.03E+6 $\pm$  3.07E+7 & 2.74E-1 $\pm$  3.79E-2\\
HNN-2-100-gelu &   4.78E-3 $\pm$ 1.14E-3  &   9.38E-2 $\pm$ 1.60E-1 &   1.47E-1 $\pm$  6.84E-2 & 7.25E-1 $\pm$  2.65E-1\\
HNN-1-200-gelu &   4.83E-3 $\pm$ 1.17E-3  &   1.26E-1 $\pm$ 2.28E-1 &  3.33E-1 $\pm$  1.18E-1 & 1.68E+0 $\pm$  3.58E-1\\
HNN-1-100-gelu &   \uwave{8.08E-4 $\pm$ 2.60E-4}  &   1.55E-2 $\pm$ 1.94E-2 &    \uwave{1.92E-2 $\pm$  1.05E-2} &  \uwave{7.38E-2 $\pm$ 3.26E-2}\\
  \colrule
  \colrule
 $\alpha$-SGHN(ours) &\textbf{4.71E-9} $\pm$ \textbf{4.26E-9}  &\textbf{1.97E-7} $\pm$ \textbf{3.57E-7}  & \textbf{1.32E-6} $\pm$ \textbf{4.05E-7}& \textbf{4.25E-6  $\pm$ 8.99E-7}\\
\end{tabular}
\end{ruledtabular}
\end{table*}

Table \ref{Toda3_com} shows the prediction comparison results of the Toda system with $N=3$. Due to the fact that there are only three particles 
(in the setting of this table) and the interactions between particles are simple, the results of MLP and HNN training are not as poor as those of the FK and Rotator systems with $N=32$, and are indeed only slightly inferior to  $\alpha$-SGHN (e.g., lower
only by 1 or 2 orders of magnitude with respect
to $\alpha$-SGHN as concerns the conservation
laws). Recall here that we avoided  carefully tuning the parameters of  $\alpha$-SGHN, which may 
indeed not be the ones leading to optimal 
performance, yet it still outperforms baseline models in several metrics. In addition, compared to the FK and Rotator models, when the number of particles is small, the width of the baseline model is correspondingly reduced to achieve the best effect. And an increase in the number of particles will greatly weaken the performance of the baseline model. However, based on graph models, the number of particles has little effect on the network model.
This strongly suggests the scalability of the
current approach to larger scale conservative
system and its comparative advantage in comparison
to earlier methodologies.
\begin{table*}
\caption{\label{Toda3_com}%
The test results for Toda system with $N=3$.  The best results are emphasized by bold fonts.    
}
\begin{ruledtabular}
\begin{tabular}{ccccc}
\textrm{}&
\textrm{Test loss}&
\textrm{Trajectory$_{MSE}$}&
\textrm{Energy$_{MSE}$}&\\
\colrule

MLP-1-25-gelu  & 2.36E-4 $\pm$ 2.93E-4 & 3.00E-3 $\pm$ 1.40E-2 & 1.56E-3 $\pm$  2.80E-3 \\
MLP-1-25-tanh  & 3.06E-7 $\pm$ 9.31E-7 & 2.14E-5 $\pm$ 8.34E-5 & 1.46E-5 $\pm$  5.97E-5 \\
MLP-1-10-silu & 1.02E-6 $\pm$ 3.58E-6 & 1.18E-5 $\pm$ 4.14E-5 & 5.57E-6 $\pm$  2.12E-5 \\
MLP-2-25-silu &  \uwave{1.58E-7 $\pm$ 4.04E-7} & 6.15E-6 $\pm$ 2.95E-5 & 2.28E-6 $\pm$  9.58E-6 \\
MLP-1-50-silu  & 4.35E-7 $\pm$ 1.50E-6 & 5.25E-6 $\pm$ 2.27E-5 & 1.02E-5 $\pm$  4.29E-5 \\
MLP-1-25-silu    & 2.53E-7 $\pm$ 4.69E-7 &  \uwave{1.28E-6 $\pm$ 1.69E-6 }&  \uwave{6.27E-7 $\pm$  2.20E-6} \\
\colrule
\colrule
HNN-1-50-gelu &  1.02E-4 $\pm$ 1.47E-4   & 2.91E-3$\pm$ 1.50E-2  &8.61E-4 $\pm$  1.58E-3 \\
HNN-1-50-tanh & 5.88E-6 $\pm$ 7.03E-6 & 6.13E-5 $\pm$ 1.83E-4  & 2.20E-6 $\pm$  4.71E-6 5\\
HNN-1-25-silu  & 8.78E-8 $\pm$ 2.42E-7 & 1.98E-6 $\pm$ 7.29E-6 & 5.82E-8 $\pm$  7.52E-8 \\
HNN-2-50-silu  & 3.04E-7 $\pm$ 1.02E-6 & 8.52E-6 $\pm$ 4.60E-5 &  \uwave{5.43E-9 $\pm$ 8.79E-9} \\
HNN-1-100-silu & 9.07E-2 $\pm$ 1.79E-1 & 6.64E-1 $\pm$ 8.55E+0 & 7.30E+1 $\pm$  3.18E+2 \\
HNN-1-50-silu &  \uwave{1.02E-7 $\pm$ 2.20E-7} &  \uwave{1.14E-6 $\pm$ 3.57E-6} & 1.31E-8 $\pm$  1.75E-8 \\
  \colrule
  \colrule
 $\alpha$-SGHN(ours)  &\textbf{7.66E-10} $\pm$ \textbf{2.53E-9}  &\textbf{ 1.25E-8} $\pm$ \textbf{2.95E-8}  & \textbf{3.64E-11} $\pm$ \textbf{3.05E-11}\\
  \hline
 \hline
 \textrm{}&
\textrm{Momentum$_{MSE}$}&
\textrm{C3$_{MSE}$}&
\textrm{}\\
\colrule

MLP-1-25-gelu   & 3.14E-3 $\pm$ 5.62E-3 & 2.04E-2 $\pm$ 3.85E-2\\
MLP-1-25-tanh  &  3.09E-6 $\pm$ 4.35E-6 &3.32E-5 $\pm$ 9.64E-5\\
MLP-1-10-silu  &  3.96E-9 $\pm$ 4.40E-9 &4.05E-6 $\pm$ 1.05E-5\\
MLP-2-25-silu  &  \uwave{5.07E-10 $\pm$ 1.03E-9} & \uwave{4.58E-7 $\pm$ 1.39E-6}\\
MLP-1-50-silu  & 1.15E-8 $\pm$ 1.65E-8 &4.61E-6 $\pm$ 1.87E-5\\
MLP-1-25-silu    &2.24E-9 $\pm$ 3.85E-9 &1.67E-6 $\pm$ 6.62E-6\\
\colrule
\colrule
HNN-1-50-gelu  &    1.96E-3 $\pm$  3.29E-3 &1.23E-2 $\pm$2.06E-2\\
HNN-1-50-tanh & 4.47E-6 $\pm$  8.55E-6 & 4.72E-5 $\pm$  8.61E-5\\
HNN-1-25-silu  & 1.45E-7 $\pm$  1.38E-7 & 4.89E-6 $\pm$  1.52E-5\\
HNN-2-50-silu   & \uwave{2.38E-8 $\pm$  4.28E-8} & 1.71E-6 $\pm$  6.57E-6\\
HNN-1-100-silu  & 1.12E+0 $\pm$  4.30E+0 & 1.11E+4 $\pm$  4.85E+4\\
HNN-1-50-silu &5.12E-8 $\pm$  7.22E-8 &  \uwave{7.96E-7 $\pm$  2.33E-6}\\
  \colrule
  \colrule
 $\alpha$-SGHN(ours) &\textbf{6.39E-11} $\pm$ \textbf{1.09E-10} & \textbf{3.96E-8} $\pm$ \textbf{1.70E-7}\\
\end{tabular}
\end{ruledtabular}
\end{table*}

\begin{table*}
\caption{\label{Toda4_com}%
The test results for Toda system with $N=4$.  The best results are emphasized by bold fonts.   
}
\begin{ruledtabular}
\begin{tabular}{cccc}
\textrm{}&
\textrm{Test loss}&
\textrm{Trajectory$_{MSE}$}&
\textrm{Energy$_{MSE}$}\\
\colrule
MLP-2-25-gelu &  1.03E-4 $\pm$ 9.57E-5   &   1.22E-3 $\pm$ 3.36E-3  &  2.39E-4 $\pm$  6.61E-4 \\
MLP-2-25-tanh  &  2.49E-5 $\pm$ 3.43E-5   &   8.08E-4 $\pm$ 3.33E-3  &  7.73E-5 $\pm$  1.83E-4 \\
MLP-2-10-silu& 7.71E-6 $\pm$ 1.69E-5   &   2.11E-4 $\pm$ 7.35E-4  &  2.07E-5 $\pm$  3.87E-5 \\
MLP-1-25-silu & 5.04E-6 $\pm$ 5.37E-6   &   3.92E-5 $\pm$ 8.15E-5  &  2.00E-5 $\pm$  3.95E-5 \\
MLP-3-25-silu& \uwave{3.80E-6 $\pm$ 6.97E-6}   &   4.28E-5 $\pm$ 1.08E-4  &  4.47E-5 $\pm$  1.69E-4 \\
MLP-2-50-silu & 5.44E-6 $\pm$ 1.11E-5 & 3.56E-5 $\pm$ 8.72E-5  & \uwave{1.31E-5 $\pm$  2.99E-5} \\
MLP-2-25-silu& 4.34E-6 $\pm$ 6.36E-6 &  \uwave{3.00E-5 $\pm$ 8.03E-5}  &4.00E-5  $\pm$ 1.33E-4 \\
\colrule
HNN-1-25-gelu &  4.56E-6 $\pm$ 1.05E-5 & 8.85E-5 $\pm$ 3.74E-4  &1.15E-6 $\pm$  1.41E-6 \\
HNN-1-25-tanh  &  2.32E-5 $\pm$ 7.00E-5 & 6.83E-4 $\pm$ 3.38E-3  &3.79E-6 $\pm$  1.28E-5 \\
HNN-1-10-silu& 3.95E-6 $\pm$ 8.89E-6 & 3.68E-5 $\pm$ 1.02E-4  &2.70E-6 $\pm$  3.43E-6 \\
HNN-2-25-silu & 3.86E-6 $\pm$ 1.42E-5 & 6.23E-5 $\pm$ 3.15E-4  &5.00E-7 $\pm$  1.86E-6 \\
HNN-1-50-silu & 2.91E-6 $\pm$ 9.00E-6 & 1.76E-5 $\pm$ 8.77E-5  &7.46E-7 $\pm$  2.66E-6 \\
HNN-1-25-silu &  \uwave{1.17E-6 $\pm$ 3.36E-6} & \uwave{8.98E-6 $\pm$ 2.66E-5} &  \uwave{1.95E-7 $\pm$ 4.22E-7} \\
\colrule
 $\alpha$-SGHN(ours) &\textbf{1.96E-9} $\pm$ \textbf{6.09E-9}  &\textbf{3.41E-8} $\pm$ \textbf{1.23E-7}  & \textbf{1.03E-8} $\pm$ \textbf{1.81E-8}\\
 \hline
 \hline
 \textrm{}&
\textrm{Momentum$_{MSE}$}&
\textrm{C3$_{MSE}$}&
\textrm{C4$_{MSE}$}\\
\colrule
MLP-2-25-gelu &  4.17E-5 $\pm$ 8.30E-5   &   6.78E-4 $\pm$ 7.61E-4  &  6.82E-2 $\pm$  9.65E-2 \\
MLP-2-25-tanh &  3.40E-5 $\pm$ 7.44E-5   &   2.14E-4 $\pm$ 4.07E-4  &  2.06E-2 $\pm$  4.33E-2 \\
MLP-2-10-silu & 7.29E-10 $\pm$ 1.83E-10   &   5.93E-5 $\pm$ 9.14E-5  &  1.09E-2 $\pm$  2.20E-2 \\
MLP-1-25-silu & 6.88E-9 $\pm$ 5.02E-9   &   9.75E-5 $\pm$ 2.27E-4  &  1.38E-2 $\pm$  3.31E-2 \\
MLP-3-25-silu & \textbf{1.15E-10 $\pm$ 1.48E-10} & 8.02E-5 $\pm$ 3.19E-4 &  1.94E-2 $\pm$  7.53E-2 \\
MLP-2-50-silu& 1.02E-9 $\pm$ 1.08E-9 & 1.57E-4 $\pm$ 3.93E-4  &1.11E-2 $\pm$  2.91E-2 \\
MLP-2-25-silu & 4.59E-10 $\pm$ 5.68E-10 &  \uwave{5.50E-5 $\pm$ 9.84E-5}  &   \uwave{8.03E-3 $\pm$ 1.80E-2 }\\
\colrule
HNN-1-25-gelu &  1.77E-6 $\pm$ 1.11E-6 & 4.04E-5 $\pm$ 6.32E-5  &1.48E-3 $\pm$  1.72E-3 \\
HNN-1-25-tanh &  7.87E-7 $\pm$ 7.56E-7 & 6.53E-4 $\pm$ 1.91E-3  &2.59E-2 $\pm$  8.09E-2 \\
HNN-1-10-silu  & 5.85E-6 $\pm$ 4.47E-7 & 1.51E-4 $\pm$ 3.17E-4  & 4.70E-3 $\pm$  7.58E-3 \\
HNN-2-25-silu & 1.29E-7 $\pm$ 1.18E-7 & \uwave{9.15E-6 $\pm$ 2.61E-5}  &  \uwave{7.34E-4 $\pm$ 2.70E-3} \\
HNN-1-50-silu & 1.19E-7 $\pm$ 1.24E-7 & 1.64E-5 $\pm$ 4.67E-5  &8.54E-4 $\pm$  2.38E-3 \\
HNN-1-25-silu &  \uwave{6.46E-8 $\pm$ 5.70E-8} &3.17E-5 $\pm$ 1.31E-4 &  7.46E-4 $\pm$ 2.64E-3\\
\colrule
 $\alpha$-SGHN(ours)&\uwave{2.11E-8  $\pm$ 3.17E-8}  &\textbf{1.64E-7} $\pm$ \textbf{2.30E-7}  & \textbf{8.33E-6} $\pm$ \textbf{1.30E-5}\\

\end{tabular}
\end{ruledtabular}
\end{table*}

\begin{table*}
\caption{\label{Toda5_com}%
The test results for Toda system with $N=5$.  The best results are emphasized by bold fonts.  
}
\begin{ruledtabular}
\begin{tabular}{ccccc}
\textrm{}&
\textrm{Test loss}&
\textrm{Trajectory$_{MSE}$}&
\textrm{Energy$_{MSE}$}&
\textrm{Momentum$_{MSE}$}\\
\colrule
MLP-1-50-gelu  &  1.87E-5 $\pm$ 2.60E-5   &  1.44E-4 $\pm$ 3.19E-4  & 9.38E-5 $\pm$  2.75E-4 &2.42E-5 $\pm$  2.71E-5\\
MLP-1-50-tanh &  3.82E-5 $\pm$ 4.72E-5   &  8.90E-4 $\pm$ 3.68E-3  & 6.44E-5 $\pm$  7.88E-5 &6.31E-5 $\pm$  5.59E-5\\
MLP-1-25-silu  &  1.23E-5 $\pm$ 1.14E-5   &  7.41E-5 $\pm$ 1.59E-4  & \uwave{1.74E-5 $\pm$  2.06E-5} &2.32E-8 $\pm$  2.80E-8\\
MLP-2-50-silu &  9.59E-6 $\pm$ 1.13E-5   &   5.59E-5 $\pm$ 1.25E-4   & 5.23E-5 $\pm$ 1.91E-4 & \textbf{5.35E-9 $\pm$  5.95E-9}\\
MLP-1-100-silu &  3.29E-1 $\pm$ 1.41E+0   &  Nan  & Nan &  Nan \\
MLP-1-50-silu &  \uwave{6.39E-6 $\pm$ 1.60E-5}   &   \uwave{5.25E-5 $\pm$ 1.25E-4}  & 2.69E-5 $\pm$ 4.05E-5 & 1.09E-8 $\pm$  8.26E-9\\
\colrule
HNN-1-25-gelu&  1.81E-4 $\pm$ 2.00E-4   & 2.20E-3 $\pm$ 7.00E-3  & 1.73E-3 $\pm$ 3.07E-3 & 4.16E-3 $\pm$ 5.77E-3\\
HNN-1-25-tanh &  4.17E-5 $\pm$ 8.75E-5   & 1.77E-3 $\pm$ 9.10E-3  & 1.74E-5 $\pm$ 2.85E-5 & 1.72E-5 $\pm$ 2.22E-5\\
HNN-1-10-silu & 7.26E-5 $\pm$ 1.21E-4   &  3.18E-3 $\pm$ 1.36E-2  & 2.32E-5 $\pm$ 5.05E-5 & 3.67E-6 $\pm$ 2.80E-6\\
HNN-2-25-silu  & 5.25E-6 $\pm$ 1.99E-5   & 1.26E-4 $\pm$ 5.90E-4  & \uwave{9.85E-7 $\pm$ 3.43E-6} & 4.15E-7 $\pm$ 3.16E-7\\
HNN-1-50-silu &  6.28E-6 $\pm$ 2.05E-5   &  1.06E-4 $\pm$ 5.45E-4  &  1.32E-6 $\pm$ 4.32E-6 &   \uwave{2.15E-7 $\pm$ 1.57E-7}\\
HNN-1-25-silu & \uwave{2.14E-6 $\pm$ 5.47E-6}   & \uwave{3.08E-5 $\pm$ 9.12E-5}  & 1.07E-6 $\pm$ 1.53E-6 & 1.36E-6 $\pm$ 9.82E-7\\
\colrule
 $\alpha$-SGHN(ours) &\textbf{2.98E-9} $\pm$ \textbf{5.19E-9}  &\textbf{4.14E-8} $\pm$ \textbf{1.79E-7}  & \textbf{7.53E-9} $\pm$ \textbf{1.19E-8}&  \uwave{1.25E-8 $\pm$ 1.81E-8}\\
 \hline
 \hline
 \textrm{}&
\textrm{C3$_{MSE}$}&
\textrm{C4$_{MSE}$}&
\textrm{C5$_{MSE}$}\\
\colrule
MLP-1-50-gelu &  1.93E-4 $\pm$ 1.59E-4   &  2.52E-2 $\pm$ 5.80E-2  &  1.92E-1 $\pm$  2.38E-1 &\\
MLP-1-50-tanh &  5.34E-4 $\pm$ 5.24E-4   &  2.96E-2 $\pm$ 3.24E-2  & 3.84E-1 $\pm$  3.95E-1 & \\
MLP-1-25-silu  &  5.43E-5 $\pm$ 7.75E-5   &  7.96E-3 $\pm$ 1.04E-2  & 7.48E-2 $\pm$  1.24E-1 & \\
MLP-2-50-silu &  \uwave{4.33E-5 $\pm$ 6.74E-5}   & 1.68E-2 $\pm$ 5.20E-2  & 8.22E-2 $\pm$ 1.68E-1 &  \\
MLP-1-100-silu &  Nan  &  Nan  & Nan \\
MLP-1-50-silu &  5.92E-5 $\pm$ 1.11E-4   &  \uwave{6.62E-3 $\pm$ 1.33E-2}  & \uwave{5.78E-2 $\pm$ 1.05E-1} &  \\
\colrule
HNN-1-25-gelu &  2.71E-2 $\pm$ 3.97E-2   & 1.29E+0 $\pm$ 2.32E+0  & 1.66E+1 $\pm$ 2.78E+1 &  \\
HNN-1-25-tanh &  3.90E-4 $\pm$ 5.42E-4   & 3.11E-2 $\pm$ 6.44E-2  & 4.78E-1 $\pm$ 8.33E-1 &\\
HNN-1-10-silu & 6.72E-4 $\pm$ 1.44E-3   &  4.28E-2 $\pm$ 9.68E-2  & 6.94E-1 $\pm$ 1.65E+0 &  \\
HNN-2-25-silu  & 2.74E-4 $\pm$ 1.12E-3   & 2.14E-2 $\pm$ 8.97E-2  & 2.96E-1 $\pm$ 1.21E+0 &  \\
HNN-1-50-silu & 6.32E-5 $\pm$ 1.86E-4   & 3.20E-3 $\pm$ 8.26E-3  & 6.14E-2 $\pm$ 1.73E-1 &  \\
HNN-1-25-silu & \uwave{2.97E-5 $\pm$ 4.69E-5}   &  \uwave{1.27E-3 $\pm$ 2.06E-3}  &  \uwave{1.96E-2 $\pm$ 3.07E-2} & \\
\colrule
 $\alpha$-SGHN(ours) &\textbf{1.13E-7} $\pm$ \textbf{1.64E-7}  &\textbf{5.91E-6} $\pm$ \textbf{9.53E-6}  & \textbf{7.99E-5} $\pm$ \textbf{1.22E-4}&\\

\end{tabular}
\end{ruledtabular}
\end{table*}

Table \ref{Toda4_com} shows the predicted comparison results of the Toda system with $N=4$.  $\alpha$-SGHN once again outperforms the baseline models in terms of test loss, and trajectory prediction. This Toda system has four conserved quantities, and, as can
be seen, the baseline model's predictions for the third and fourth conserved quantities are 
progressively getting worse, while $\alpha$-SGHN still performs better. This may be because SGHN 
retrieves the underlying structure of the system, making predictions for various quantities more stable. In an interesting outlier, in comparison
to what was seen in earlier tables, for the momentum prediction, the MLP performs the best, but the corresponding result of $\alpha$-SGHN is also sufficient and not too far worse. Overall,
however, the comparative advantage of
$\alpha$-SGHN when one views all the relevant
diagnostics is demonstrable.

Table \ref{Toda5_com} shows the predicted comparison results of the Toda system with $N=5$. The system contains five conserved quantities, and it is evident that the baseline model cannot guarantee the conservation of the relevant quantities. This
issue is especially prominent for the third, fourth, and fifth conserved quantities. Indeed, here we
see a degradation of MLP and HNN the higher
the relevant conservation law. 
$\alpha$-SGHN predicts the interactions between system particles, making the learning of conserved quantities more stable.  By the way, once again,
the momentum conservation law is still not optimal
in the case of the $\alpha$-SGHN, yet the overall
evidence is even more overwhelmingly in favor
of our discovery of the network structure.
 
\begin{figure*}%
\centering
\includegraphics[scale=0.63]{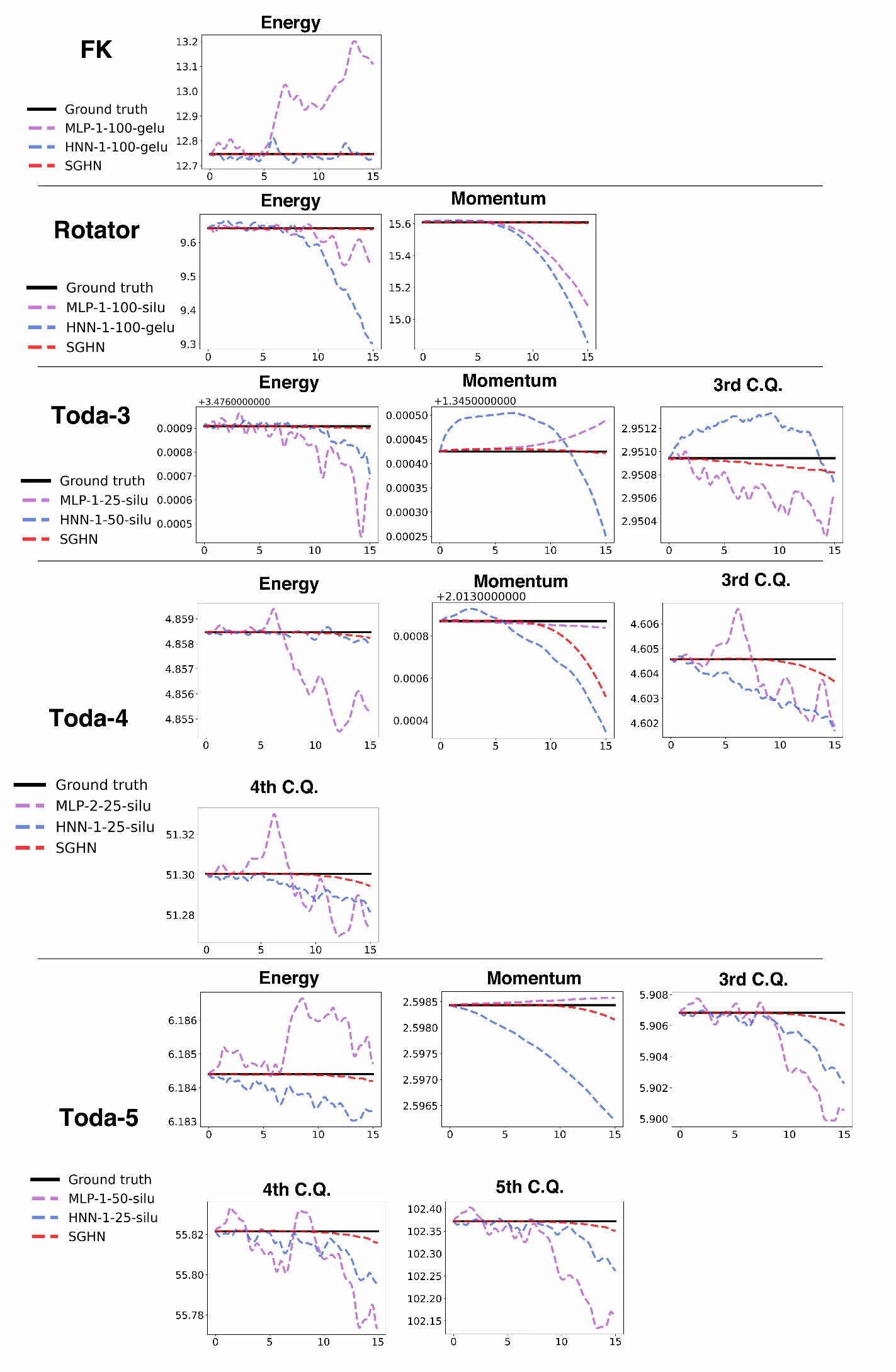}
\caption{The evolution over time of the average true conservation values and the average predicted conservation values of 20 samples. C.Q. represents conserved quantity.}\label{c-frt}
\end{figure*}

\begin{figure*}%
\centering
\includegraphics[scale=0.8]{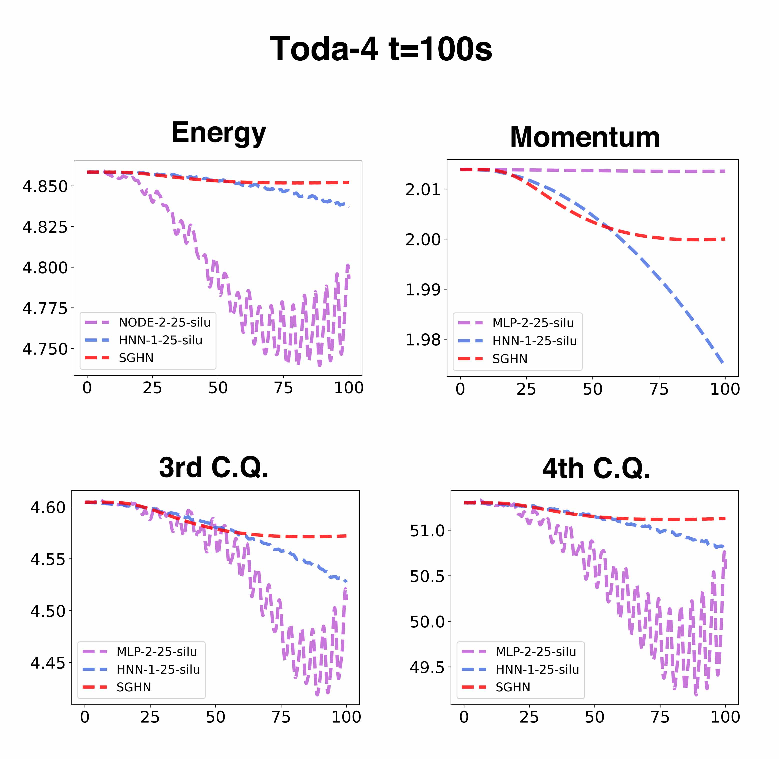}
\caption{The evolution over time of the average true conservation values and the average predicted conservation values of 20 samples. C.Q. represents conserved quantity.}\label{c2-frt}
\end{figure*}

Fig. \ref{c-frt} shows the evolution of the average true conservation law values and average predicted conservation values over time for 20 samples. 
The method of obtaining 20 samples can be found in Section III.B.
The baseline model adopts the network configuration with the best performance. Please note that the range of the $y$- axis is different. It can be seen that $\alpha$-SGHN has the best performance. The $\alpha$-SGHN prediction of momentum in Toda with $N=4$ appears to have a significant attenuation, but its error variation is around 1e-4.
Generally, over the time scales of a few 10s of units considered, it is clear that the $\alpha$-SGHN has the weakest fluctuations. 
Fig. \ref{c2-frt} shows the evolution of Toda-4's predicted conservation laws over time when the prediction time increases by 20 times the training time. It can be seen that overall, $\alpha$-SGHN performs better than the baseline model (when taking into account all the conserved quantities and their
deviation from the original value).

\subsection{Why graph neural networks have advantages over traditional neural networks}

\begin{figure*}%
\centering
\includegraphics[scale=0.7]{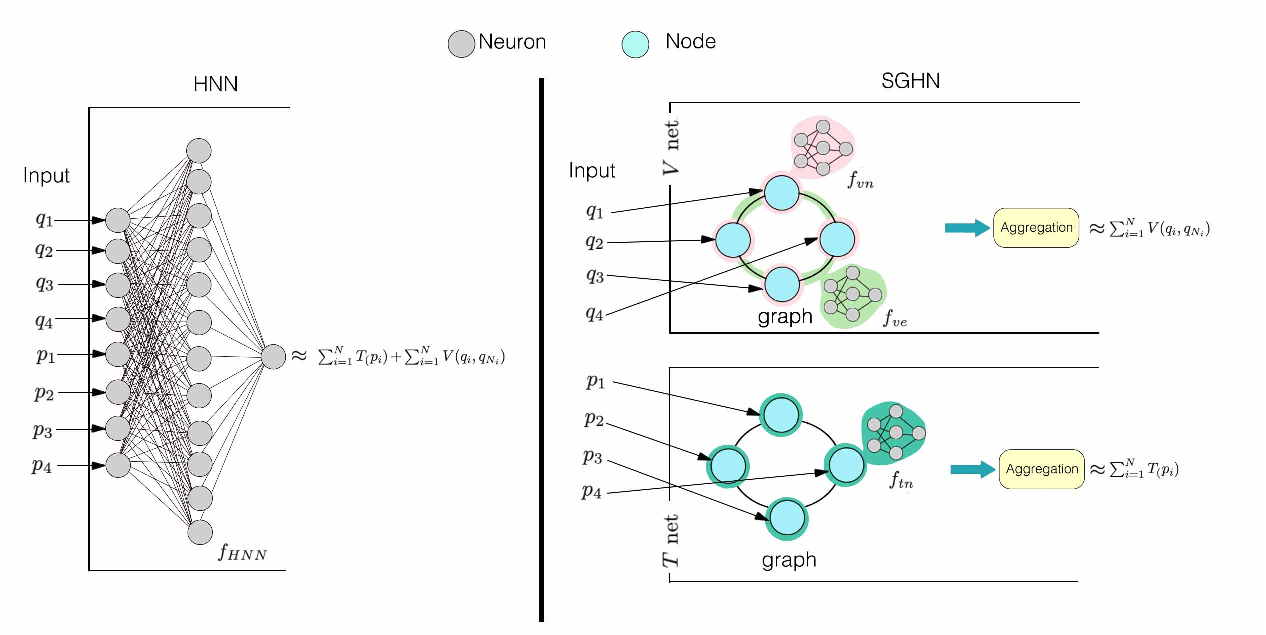}
\caption{Comparison of network structures for learning energy part between HNN and SGHN.}\label{fig_neural}
\end{figure*}

In this section, we provide an intuitive explanation of why graph neural network-based methods are superior to traditional neural networks. Fig. \ref{fig_neural} shows a comparison of network structures for learning energy part between HNN and SGHN.

On the left is a HNN, which is mainly learned by an MLP denoted as $f_{HNN}$ to learn the Hamiltonian of the system,   
 That is to say, it uses $f_{HNN}(q_1,\cdots,q_N,p_1,\cdots,p_N)$ to approximate $H=\sum_{i=1}^N T{(p_i)}+\sum_{i=1}^N V({q_i,q_{N_i}})$, where the input is $q_1,\cdots,q_N,p_1,\cdots,p_N$. From the  Fig. \ref{fig_neural} left, it can be seen that there is no connection between the input data and they are connected to the neurons of the next layer of the network in the same way.

 On the right is the  $\alpha$-SGHN part 2. It can be seen that the input data are first assigned to the system graph structure learned in the first part (the graph structure encodes the interaction relationships between particles in the system). Then, the information of the nodes and edges is aggregated through MLPs placed on the nodes and edges. That is to say, $\alpha$-SGHN knows more about the interaction information of the system than HNN, hence it is more effective than HNN.

\subsection{Comparison with Graph Neural Networks}

In this section, we will compare our method with two graph neural network methods Hamiltonian ODE graph network (HOGN)\cite{sanchez2019hamiltonian}, Hamiltonian graph neural networks (HGNN)\cite{bishnoi2023learning}, and the parameter settings of the graph neural network model are the same as those in~\cite{bishnoi2023learning}. 
We define the initial conditions as follows:
 \begin{eqnarray}
&&\textbf{q}_0(i)=\lambda_{i} \sin\bigg(\frac{(i-1) \pi}{N-1}\bigg),\\
&&\textbf{p}_0(i)=0,\quad i=1,\cdots,N,
\end{eqnarray}
where $\lambda_{i} \sim \mathcal{U}(0,1)$. 
Under the same input conditions,i.e. without knowing the links of the particles, a fully connected graph is used like $\alpha$-SGHN (assuming that particles all have interactions between them), and the test results are shown in Table \ref{com_ghnn}.
It can be seen that other graph based methods cannot work without knowing the particle links within the system (i.e., the particle interactions). Our method works because it can predict the particle links of the system,
uncovering the matrix of the inter-particle interactions.
In the case of known lattice particles links,
the accuracy of  HOGN and HGNN cannot support to 
 learn systems with non-even symmetric potential energy, see Table \ref{com_ghnn2}. 
\begin{table*}
\caption{\label{com_ghnn}%
The test results for HGNN, HOGN and $\alpha$-SGHN. The best results are emphasized by bold fonts. The link relationship of the lattice system is unknown, so a fully connected network is used.
}
\begin{ruledtabular}
\begin{tabular}{ccccc}
\multicolumn{5}{c}{\textrm{FK system with $N=32$}}\\
\textrm{}&
\textrm{Test loss}&
\textrm{Trajectory$_{MSE}$}&
\textrm{Energy$_{MSE}$}&
\textrm{Momentum$_{MSE}$}\\
\colrule
HOGN  & 3.35E-2 $\pm$ 5.82E-2  &     1.26E+0 $\pm$ 4.59E+0 &   4.77E+3 $\pm$ 1.43E+4 & -- \\
HGNN  & 3.44E-2 $\pm$ 6.44E-3  &     1.82E-1 $\pm$ 9.33E-2  &   8.60E-2 $\pm$ 1.06E-1 & -- \\
 $\alpha$-SGHN   &\textbf{2.46E-9} $\pm$ \textbf{1.61E-9}  &\textbf{2.73E-8} $\pm$ \textbf{3.67E-8}  & \textbf{1.72E-8} $\pm$ \textbf{2.14E-8}& --\\
  \hline
 \hline
 \multicolumn{5}{c}{\textrm{Rotator system with $N=32$}}\\
 \textrm{}&
\textrm{Test loss}&
\textrm{Trajectory$_{MSE}$}&
\textrm{Energy$_{MSE}$}&
\textrm{Momentum$_{MSE}$}\\
\colrule
HOGN  & 1.90E-1 $\pm$ 5.23E-2  &    Nan &    Nan  & Nan   \\
HGNN  & 1.43E-1 $\pm$ 3.04E-2  &     3.02E-1 $\pm$ 1.21E-1 & 8.92E-1 $\pm$ 8.31E-1 & 2.38E+1 $\pm$ 7.93E+0 \\
 $\alpha$-SGHN &\textbf{3.90E-9} $\pm$ \textbf{1.66E-9}  &\textbf{3.93E-8} $\pm$ \textbf{5.08E-8}  & \textbf{5.27E-8} $\pm$ \textbf{2.13E-8}& \textbf{8.07E-8} $\pm$ \textbf{3.03E-8}\\
   \hline
 \hline
 \multicolumn{5}{c}{\textrm{Toda system with $N=5$}}\\
 \textrm{}&
\textrm{Test loss}&
\textrm{Trajectory$_{MSE}$}&
\textrm{Energy$_{MSE}$}&
\textrm{Momentum$_{MSE}$}\\
\colrule
HOGN  & 3.44E-2 $\pm$ 3.11E-2  &     1.89E-1 $\pm$ 1.94E-1 &   2.92E-2 $\pm$ 5.20E-2  &   1.92E-1 $\pm$ 2.62E-1 \\
HGNN  & 1.37E-2 $\pm$ 1.01E-2  & 2.00E-1 $\pm$ 1.84E-1 &  4.81E-3 $\pm$ 5.92E-3&  6.17E-3 $\pm$ 8.51E-3  \\
 $\alpha$-SGHN  &\textbf{4.80E-9} $\pm$ \textbf{6.40E-9}  &\textbf{8.61E-8} $\pm$ \textbf{1.10E-7}  & \textbf{7.56E-10} $\pm$ \textbf{6.83E-10}  & \textbf{1.06E-9} $\pm$ \textbf{1.50E-9}\\
\end{tabular}
\end{ruledtabular}
\end{table*}

\begin{table*}
\caption{\label{com_ghnn2}%
In the case of known lattice system links, HOGN and HNN cannot learn systems with non even symmetric potential energy.
}
\begin{ruledtabular}
\begin{tabular}{ccccc}
  \multicolumn{5}{c}{\textrm{Toda system with $N=5$}}\\
 \textrm{}&
\textrm{Test loss}&
\textrm{Trajectory$_{MSE}$}&
\textrm{Energy$_{MSE}$}&
\textrm{Momentum$_{MSE}$}\\
\colrule
 
HOGN  with graph & 2.42E-3 $\pm$ 2.66E-3  &     1.20E-2 $\pm$ 1.62E-2 & 2.61E-3 $\pm$ 3.78E-3 & 2.67E-3 $\pm$ 2.45E-3 \\
HGNN with graph &2.62E-3 $\pm$ 2.93E-3  &     5.30E-3 $\pm$ 7.25E-3 &  1.97E-3 $\pm$3.35E-3 & 1.88E-4 $\pm$ 2.90E-4  \\
  $\alpha$-SGHN &\textbf{3.90E-9} $\pm$ \textbf{1.66E-9}  &\textbf{3.93E-8} $\pm$ \textbf{5.08E-8}  & \textbf{5.27E-8} $\pm$ \textbf{2.13E-8}& \textbf{8.07E-8} $\pm$ \textbf{3.03E-8}\\
 \end{tabular}
\end{ruledtabular}
\end{table*}

{{ 
\section{The relationship between the complete integrability of a system and the learning performance of neural networks  }
In this section, we explore the impact of the system's conservation laws on the predictive performance of neural network models. We construct a hybrid system of FK and Toda (FK-Toda), as shown below: 
 \begin{eqnarray}\label{phieq2}
H&=&\sum_{i=1}^N\frac{p_i^2}{2}+\mu\sum_{i=1}^N  \exp(q_i-q_{i+1}) \nonumber\\
&+&(1-\mu)\Bigg( \frac{(q_{i+1}-q_i)^2}{2}+1-\cos(q_i)\Bigg).
\end{eqnarray}
When $\mu=0$, it is an FK system with a conserved quantity, which is energy. When $\mu=1$, it is a Toda system, completely integrable, with $N$ conserved quantities. 
Notice that this is an interesting dynamical system in its own right, inspired by the well-known Salerno model in the context of discrete nonlinear Schr{\"o}dinger equations~\cite{salerno1992quantum}.
The latter has been explored recently in the vicinity of the integrable limit as concerns the sensitivity of higher-order conserved quantities~\cite{mithuno}
and the ``detectability'' of a model's integrable nature.
Here, we explore the predictive performance of the FK-Toda system and how it changes for the model when $\mu \in [0,1]$. 

For evaluation, we use the Mean Absolute Percentage Error (MAPE) which is defined as
 \begin{eqnarray}\label{MAPE}
 MAPE=\frac{1}{N_T} \sum_{t=1}^{N_T} \frac{1}{N} \sum_{i=1}^N\left|\frac{\hat{\mathbf{F}}_{t, i}-\mathbf{F}_{t, i}}{\mathbf{F}_{t, i}}\right|,
\end{eqnarray}
where $\hat{\mathbf{F}}_{t, i}$ and $\mathbf{F}_{t, i}$ are the predicted value and the corresponding ground truth at the $t$-th time for the $i$-th particle.

Fig. \ref{fig_todafk1} shows the changes in the predicted MAPE of trajectory and energy of the network model as $\mu$ increases when $N=5$. 
In order to display this clearly, $y$-axis using semilog scale.
It can be seen that the MAPE of $\alpha$-SGHN has always been the lowest.
As $\mu$ increases, the system becomes more complex and the predictive performance may decrease. However, when  $\mu$ approaches 1, the predictive performance of the model improves.  when $\mu=1$, we retrieve the completely integrable model, and there we can see how the $\alpha$-SGHN predictions visibly capture the integrable structure; 
We observed similar MAPE trends for HNN and MLP. That is to say, when complex systems approach integrability, the predictive performance of the model will improve, and when the system is completely integrable, the predictive performance of the model will reach its best. This is an interesting
feature worth exploring further in the future (e.g., its
potential validity in other settings).


\begin{figure*}%
\centering
\includegraphics[scale=0.5]{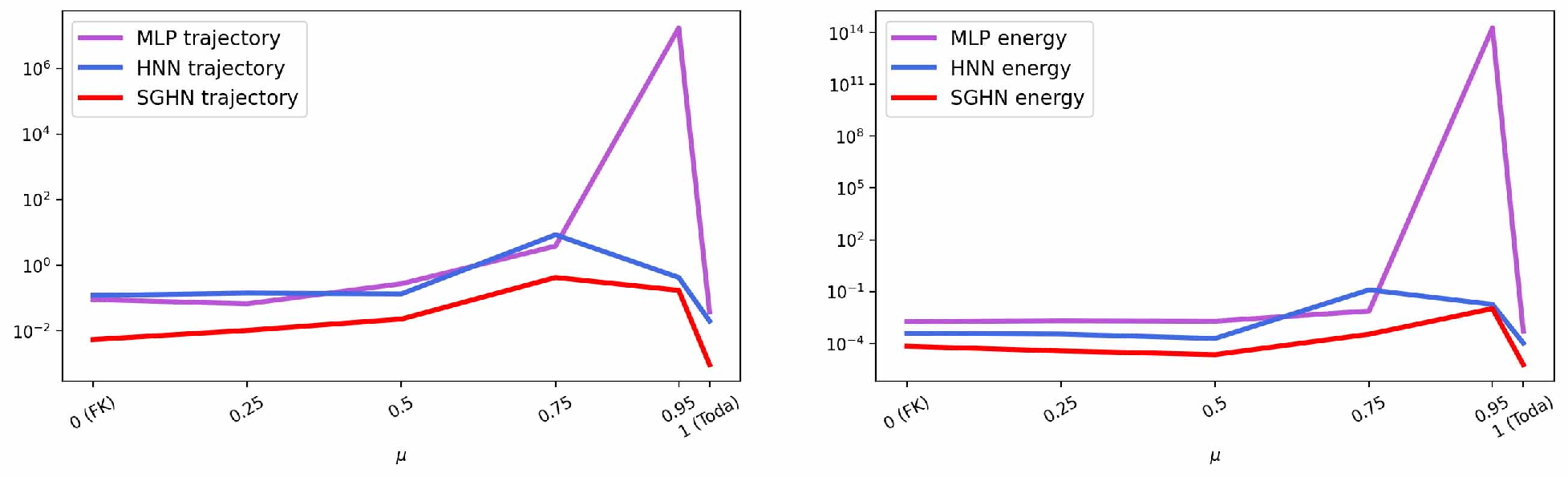}
\caption{The mean of MAPE of 20 test samples when $\mu$ changes from 0 to 1 for FK-Toda system with $N=5$. HNN has a hidden layer with 25 units per layer, while MLP has a hidden layers with 50 units per layer. The activation functions are all silu. In order to display this clearly, $y$-axis using semilog scale. }\label{fig_todafk1}
\end{figure*} 

To avoid the above conclusion being an exception, we also constructed in the same spirit a hybrid system of the Fermi–Pasta–Ulam-Tsingou \cite{fermi1955studies,berman2005fermi,VAINCHTEIN2022133252} and Toda (FPUT-Toda):
 \begin{eqnarray}\label{phieq1}
H&=&\sum_{i=1}^N\frac{p_i^2}{2}+\mu\sum_{i=1}^N  \exp(q_i-q_{i+1}) \nonumber\\
&+&(1-\mu)\Bigg( \frac{(q_{i+1}-q_i)^2}{2}+\frac{(q_{i+1}-q_i)^3}{6}\Bigg).
\end{eqnarray}
When $\mu=0$, it is an FPUT system with two conserved quantities, which are energy and momentum. When $\mu=1$, it is a Toda system, completely integrable, with $N$ conserved quantities.  

Fig. \ref{fig_todafpu1} shows the results of the FPUT-Toda system testing.  The MAPE of $\alpha$-SGHN has always been the lowest. 
When the system is fully integrable, i.e. $\mu=1$, the performance of the three neural network models becomes optimal.
Overall, our results indicate that neural networks have better performance in modeling fully integrable systems, as also discussed above.
 
\begin{figure*}%
\centering
\includegraphics[scale=0.5]{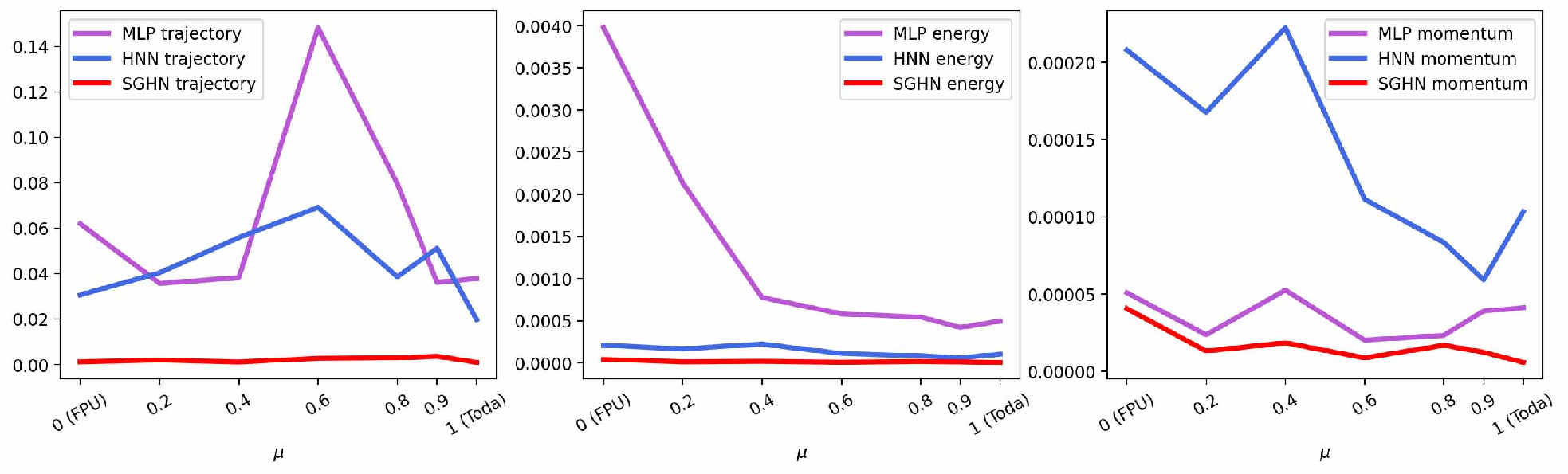}
\caption{The mean of MAPE of 20 test samples when $\mu$ changes from 0 to 1 for FPU-Toda system with $N=5$. HNN has a hidden layer with 25 units per layer, while MLP has a hidden layers with 50 units per layer. The activation functions are all silu.}\label{fig_todafpu1}
\end{figure*} 

\section{Conclusions and Future Challenges}
In this work, we introduced a new model $\alpha$-SGHN  for studying lattice systems. Firstly, it can capture the interactions between particles in complex, strongly nonlinear Hamiltonian lattice systems with multiple degrees of freedom. Secondly, $\alpha$-SGHN can utilize the predicted particle link for further and more accurate system behavior prediction. It eliminates the limitation of knowing the system structure in advance when studying lattice systems based on graph neural networks, and provides a new method for inferring particle interactions based solely on particle trajectories.
$\alpha$-SGHN can discover particle interactions in lattice systems solely from their motion trajectories, without requiring any prior knowledge of structural information.
In addition, we also investigated whether the particle behavior predicted by $\alpha$-SGHN  preserves the system's conservation laws. We tested a variety of systems ranging from one and two conserved quantities, to controllably many conserved quantities and ultimately to (in the case of as many as the 
degrees of freedom, the scenario of) complete integrability. The experimental results show that the trajectory predicted by $\alpha$-SGHN preserves all the conserved variables of the system, and generically features a performance that is far superior to the baseline model, although we did find isolated exceptions of other networks performing better (e.g., the momentum conservation in the some of the Toda lattice examples). 
Nevertheless, overall, the performance
of $\alpha$-SGHN  was found to be considerably
superior to that of the baseline models to which
it was compared herein, lending considerable
promise to it for future applications.
This may be due to the graph structure in $\alpha$-SGHN implying information about particle interactions.

These findings seem to suggest the relevance 
of extending considerations of the $\alpha$-SGHN
approach to a variety of other models, including
ones with beyond nearest-neighbor interactions
and even ones associated with long-range
interactions and the related fractional
derivatives, a number of examples of which have been
recently explored in~\cite{cuevas}.
This may be quite useful in not only discovering
such models and reconstructing their trajectories,
but also towards examining their potential
conservation laws and exploring their integrability, 
as was
done herein in a range of benchmark cases.
Such studies are currently in progress and will be
reported in future publications.

\section*{Acknowledgments}
 This material is based upon work supported by the
U.S. National Science Foundation under the awards
PHY-2110030, DMS-2204702 (PGK), and by
 National Science Foundation of China under the awards Grant Nos. 12371187, 11871140.


\appendix
\section{Supplementary experiment on structural prediction}
To explore the universality of $\alpha$-SGHN, we conducted data analysis on systems with long-range interactions in two dimensions to see if particle interactions can be accurately learn.

We take the Klein-Gordon (KG) lattice system with long-range interactions (KG-LRI) \cite{koukouloyannis2013multibreathers}  and the 2D FK system  \cite{granato1999dynamical}. The methods for obtaining datasets and predicting trajectory results have been extensively studied in previous work, see \cite{geng2024separable}. Here we will only explore the results of structural learning.

The KG-LRI system is described by
\begin{eqnarray}
H&=&\sum_{i=1}^N \Bigg(\frac{p_{i}^2}{2}+ a\frac{(q_{i+1}-q_{i})^2}{2}+b\frac{(q_{i+2}-q_{i})^2}{2}\nonumber \\ 
&+&\frac{q_{i}^2}{2}+\frac{q_{i}^4}{4}\Bigg),	
\end{eqnarray}
where $a>0$, $b\geq 0$. We extend our investigation to the periodic KG-LRI model, presuming $a=b=1$ and a total of $N=32$ nodes.

The 2D FK system is described by
\begin{eqnarray}
H&=&\sum_{i,j=1}^N \Bigg(\frac{p_{i,j}^2}{2}+ a\frac{(q_{i+1,j}-q_{i,j}-\rho)^2}{2}+b\frac{(q_{i,j+1}-q_{i,j})^2}{2}\nonumber \\ 
&-&\cos(q_{i,j})\Bigg),	
\end{eqnarray}
 where $a>0$, $b>0$, and $\rho$ denotes the average particle distance absent an external potential.
Within our 2D FK model framework, we set $a=b=\rho=1$. Assuming $M=N=12$, this model portrays a quadrilateral grid with a total of 144 particles.

In fact, in our model, the link is easily learned. We set epoch=1200, and the learning rate drops to 1e-5 after 4600 epochs. The other settings of the model remain unchanged, the same as before.

Fig. \ref{fig_edges_app} shows $|\alpha|$ and the particle interaction relationship inferred from  $|\alpha|$. We have displayed partial particle IDs in the 2D FK system. Obviously, these systems satisfy the periodic boundary conditions.
 
\begin{figure*}%
\centering
\includegraphics[scale=0.63]{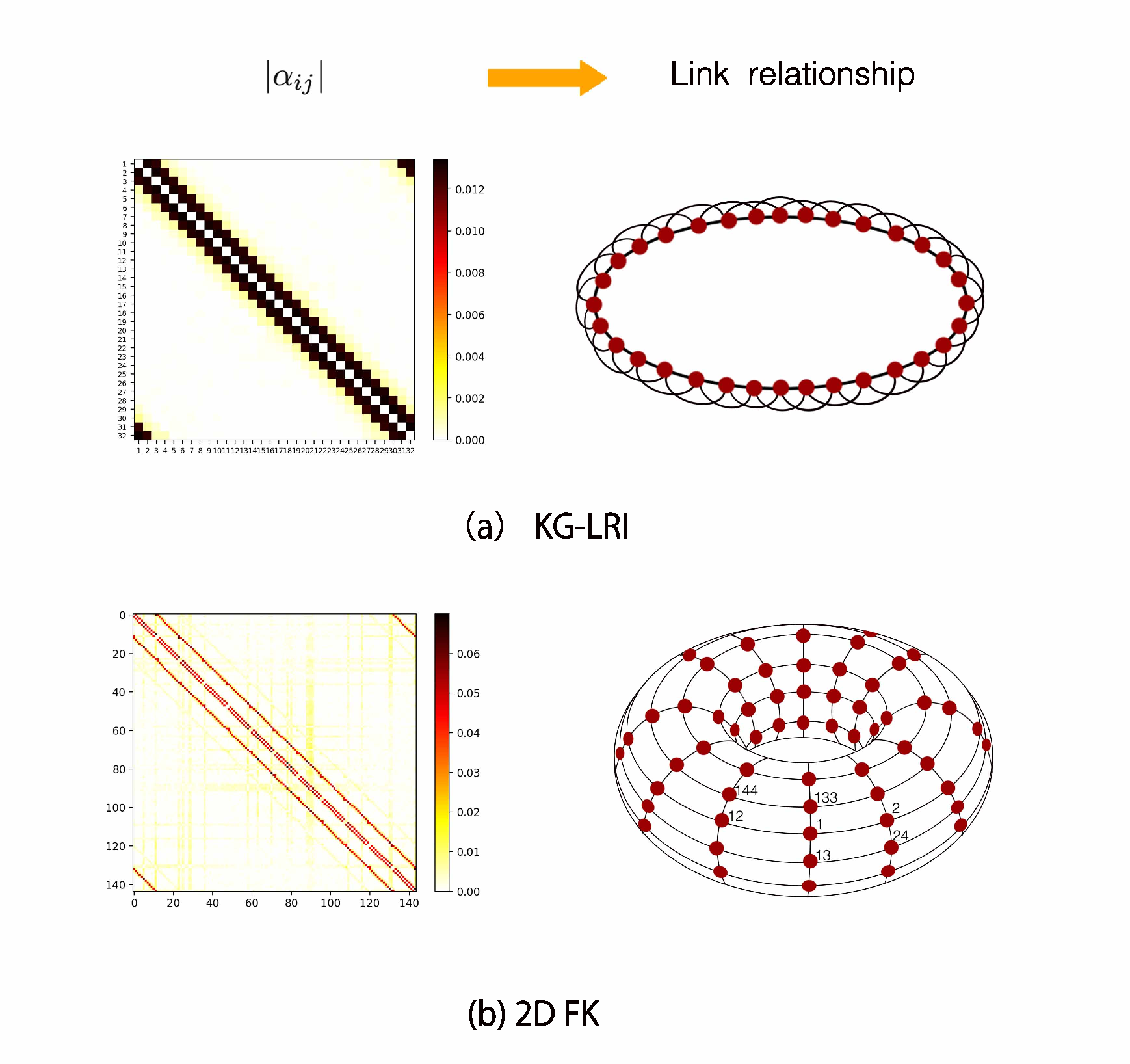}
\caption{Left: the $x$-axis and $y$-axis represent nodes, and the color coded $|\alpha_{i,j}|$ values. Right: The particle link relationship is extracted from the left. It satisfies the periodic boundary conditions.  }\label{fig_edges_app}
\end{figure*} 
\normalem
\bibliography{apssamp}

\end{document}